\documentclass[11pt]{amsart}
\usepackage[T1]{fontenc}
\usepackage{amsmath,amssymb,mathtools}
\usepackage{needspace}
\usepackage[unicode,hidelinks]{hyperref}
\numberwithin{equation}{section}
\newtheorem{theorem}{Theorem}[section]
\newtheorem{lemma}[theorem]{Lemma}
\newtheorem{proposition}[theorem]{Proposition}
\newtheorem{corollary}[theorem]{Corollary}
\theoremstyle{remark}

\newcommand{\F}{\mathbb F}

\newcommand{\R}{\mathbb R}
\newcommand{\C}{\mathbb C}
\newcommand{\PP}{\mathbb P}
\DeclareMathOperator{\MaxR}{MaxR}

\DeclareMathOperator{\im}{im}
\title[Maximal tensor rank]{The asymptotics of 
	maximal rank of three-way tensors}
\author{Jaros{\l}aw Buczy{\'n}ski}
\address{Jarosław Buczyński\\
Institute of Mathematics of the
  Polish Academy of Sciences\\
  ul. Śniadeckich 8\\
  00-656 Warszawa, Poland
}
\email{jabuczyn@impan.pl}

\author{Mateusz Micha{\l}ek}
\address{Mateusz Micha{\l}ek\\
University of Konstanz\\
78457 Konstanz, Germany\\
}
\email{mateusz.michalek@uni-konstanz.de}

\date{22 November 2026}
\subjclass[2020]{Primary 15A69; Secondary 14N07, 14M17}
\keywords{Tensor rank, maximal rank, Segre varieties, regular decompositions of zero}
\usepackage[textsize=tiny]{todonotes}
\usepackage{lmodern}
\usepackage{microtype}
\microtypesetup{expansion=false}

\begin{document}
\begin{abstract}
We prove that the maximal rank of an $n\times n\times n$ tensor is
$n^2/3+O(n^{3/2})$, over the real numbers and over every algebraically
closed field. More generally, in format
$\lfloor an\rfloor\times\lfloor bn\rfloor\times\lfloor cn\rfloor$,
with $a,b,c>0$ fixed, the leading term is $\frac{abc}{a+b+c}n^2$.
The argument starts with a criterion for rank with respect to an
arbitrary projective variety: a decomposition of zero with surjective
differential, from a smooth parameter space, bounds the rank of every
point. For Segre varieties, we construct such decompositions using a
basis of units in a split algebra. We also give an optimized finite bound and a closed estimate with explicit constants for the maximal rank.
\end{abstract}
\maketitle
\section{Introduction}

How much larger can the rank of an exceptional tensor be than the rank
of a general tensor? For matrices there is no difference. For tensors
of order three, the two ranks can differ, and a generic decomposition
does not give a decomposition of every tensor. We show that this
difference disappears in the leading term when the dimensions grow
in fixed proportions.

Let $\F$ be a field. For $T\in\F^{n_1}\otimes\F^{n_2}\otimes\F^{n_3}$, let $R_\F(T)$ be
the least number of decomposable summands in an expression for $T$.
We write
\[
 \MaxR_\F(n_1,n_2,n_3)
 =\max \{R_\F(T): T\in\F^{n_1}\otimes\F^{n_2}\otimes\F^{n_3}\}.
\]

\Needspace{9\baselineskip}
\begin{theorem}\label{thm:segre}
Suppose $\F=\R$, or $\F$ is an algebraically closed field. For fixed real numbers $a,b,c>0$,
\begin{equation}\label{eq:segre-main}
 \MaxR_\F(\lfloor an\rfloor,\lfloor bn\rfloor,\lfloor cn\rfloor)
 =\frac{abc}{a+b+c}n^2+O_{a,b,c}(n^{3/2}).
\end{equation}
In particular, $\MaxR_\F(n,n,n)=n^2/3+O(n^{3/2})$.
\end{theorem}
In Lemma \ref{lem:effective} we give explicit estimates of the error term $O(n^{3/2})$.

Another interpretation of our results compares rank and border rank. We prove that for tensors in $(\F^n)^{\otimes 3}$ the ratio of maximal rank to maximal border rank tends to one as $n$ tends to infinity. It is known that for special tensors of arbitrary large format the ratio of rank to border rank can exceed two \cite{LandsbergMichalekAbelian, ZuiddamGap}.

\subsection*{Earlier bounds}
The leading coefficient in \eqref{eq:segre-main} is suggested by
parameter counting. Strassen and Lickteig established the corresponding
generic-rank estimates \cite{Lickteig,Strassen}. The distinction between
rank and border rank is discussed in
\cite{BuczynskiLandsberg,LandsbergBook,LandsbergMichalek}.

Blekherman and Teitler prove that maximal rank with respect to a
projective variety is at most twice generic rank
\cite{BlekhermanTeitler}.
Over an algebraically closed field of characteristic zero,
Buczy{\'n}ski, Han, Mella, and Teitler improve this to $2g-1$ for a
projective homogeneous variety of generic rank $g$, and to $2g-2$
when its last nontrivial secant variety is a hypersurface
\cite[Thm~3.9]{BuczynskiHanMellaTeitler}.
For three-way tensors, the classical estimates of Atkinson--Stephens
and Atkinson--Lloyd, and their real counterparts due to
Sumi--Miyazaki--Sakata, give the stronger bound
\begin{equation}\label{eq:classical}
 \MaxR_\F(n,n,n)\le\frac{n(n+1)}2,
 \qquad \F=\R\text{ or }\C;
\end{equation}
see \cite{AtkinsonLloyd,AtkinsonStephens} and
\cite[Thms~3.4 and 3.7]{SumiMiyazakiSakata}.
Thus Theorem~\ref{thm:segre} replaces the coefficient $1/2$ in this
uniform upper bound by the optimal parameter-count coefficient $1/3$.

\subsection*{The argument}
The proof is organized around a local criterion with a global
consequence. Suppose a sum of rank-one points is zero and the
differential of the sum, taken on a smooth parameter space, is
surjective. The image then contains a neighbourhood of zero.
Every vector has a nonzero scalar multiple in this neighbourhood,
so every vector has rank at most the number of summands. We state
this observation for rank with respect to an arbitrary projective
variety in Lemma~\ref{lem:regular-zero}.

For Segre varieties, we combine a generic decomposition in one
format with a special decomposition coming from the algebra structure on $\F^k$. The tensor product decomposition we obtain turns out to have sum zero and surjective differential, exactly thanks to 
the spanning and cancellation conditions we introduce
for the smaller decomposition.
This yields Theorem~\ref{thm:segre}. Section~\ref{sec:effective}
optimizes the integer parameters and gives, for $n\ge64$,
\[
 \MaxR_\F(n,n,n)\le
 \frac{n^2}{3}+\frac{2\sqrt6}{9}n^{3/2}+\frac83n^{5/4}+\frac{65}{12}n.
\]
The tangent-space methods to prove surjectivity of the differential map, including their characteristic-independent
formulations, are developed in Appendix~\ref{app:tangents}.

\section{A regular cancellation bounds \texorpdfstring{$X$}{X}-rank}\label{sec:criterion}

Let $W$ be a finite-dimensional vector space over $\F$, and let
$X\subseteq\PP(W)$ be a projective variety with affine cone
$\widehat X\subseteq W$. We assume that $\widehat X(\F)$ linearly spans $W$.
For $w\in W$, define
\[
 R_{X,\F}(w)=\min\{r\in\mathbb Z_{\ge0}:w=x_1+\cdots+x_r,\ x_i\in\widehat X(\F)\},
\]
and write $\MaxR_\F(X)=\max_{w\in W}R_{X,\F}(w)$.
Here $R_{X,\F}(0)=0$. This is the usual projective $X$-rank of
$[w]$ when $w\ne0$. In particular, rank is unchanged by multiplying
$w$ by a nonzero scalar. Our main case of interest is the Segre variety $X=\PP(V_1)\times \PP(V_2)\times \PP(V_3)\subset \PP(V_1\otimes V_2\otimes V_3)$, where $X$-rank is the usual tensor rank.

%\Needspace{14\baselineskip}
\begin{lemma}[Regular-cancellation criterion]\label{lem:regular-zero}
Let $\F=\R$, or let $\F$ be algebraically closed. Let
$\phi:Y\to\widehat X$ be a morphism of varieties. For an integer $r\ge1$, consider
\[
 \Phi_r:Y^r\longrightarrow W,
 \qquad (y_1,\ldots,y_r)\longmapsto\sum_{i=1}^r\phi(y_i).
\]
Suppose there exists a tuple $(y_1,\dots,y_r)\in Y(\F)^r$ of smooth points of $Y$, such that
$\Phi_r(y_1,\ldots,y_r)=0$, and $\mathrm{d}\Phi_r$ is surjective at $(y_1,\dots,y_r)$. Then
\[
 \MaxR_\F(X)\le r.
\]
In particular, it suffices to find smooth points
$x_i\in\widehat X(\F)$ such that
$\sum_i x_i=0$ and $\sum_i T_{x_i}\widehat X=W$.
\end{lemma}
\begin{proof}
First let $\F$ be algebraically closed. Both $Y^r$ and $W$ are
smooth at the relevant points. The differential criterion for
smoothness shows that $\Phi_r$ is smooth near the given tuple.
Smooth morphisms are open
\cite[\href{https://stacks.math.columbia.edu/tag/056G}{Lemma~056G}]{Stacks}, %\jb{zmieniłem formę odnośnika}
so the
image contains a Zariski-open neighbourhood $O$ of zero. Every
$\F$-point of $O$ has an $\F$-point in its fibre: the fibre is a
nonempty scheme of finite type over an algebraically closed field.
Thus $R_{X,\F}(w)\le r$ for $w\in O(\F)$.

For $w\ne0$, the line $\F w$ meets $O$ in an open set containing
zero. Since $\F$ is infinite, this set contains $\lambda w$ with
$\lambda\ne0$. Hence $R_{X,\F}(w)=R_{X,\F}(\lambda w)\le r$.
Over $\R$, the submersion theorem gives a Euclidean neighbourhood
of zero in the image, and the same scaling argument applies.
For the final assertion, take $Y=\widehat X_{\mathrm{sm}}$ and
$\phi$ to be the inclusion.
\end{proof}

%We call the tuple in the lemma a \emph{regular zero} of $\Phi_r$.
Following the notation of the lemma, we call an $r$-tuple $(y_1,\dots,y_r)$ a \emph{cancellation} if it belongs to $\Phi_r^{-1}(0)$. We call an $r$-tuple \emph{regular} if the differential at $(y_1,\dots,y_r)$ is surjective.  

We note that in the lemma no smoothness of $X$ is assumed; smoothness is required only at the
chosen points of the parameter space. It is essential that the tuple
maps to zero rather than to an arbitrary point.

For the Segre variety $X=\PP(V_1)\times\PP(V_2)\times\PP(V_3)$,
use the parameter space $Y=V_1\times V_2\times V_3$ and
$\phi(v_1,v_2,v_3)=v_1\otimes v_2\otimes v_3$. At a tuple
$v=(v_{i,j})$, the differential of its $r$-term sum has image
\begin{align}
 \im(\mathrm{d}\Phi_r|_v)=\sum_{i=1}^r\bigl(
 &V_1\otimes v_{i,2}\otimes v_{i,3}
 +v_{i,1}\otimes V_2\otimes v_{i,3}+v_{i,1}\otimes v_{i,2}\otimes V_3\bigr).
 \label{eq:segre-differential}
\end{align}
%This formula will also identify the tangent-space input needed below.

\begin{corollary}[Segre case]\label{cor:segre-zero}
Let $\F=\R$, or let $\F$ be algebraically closed. If
$\sum_{i=1}^r v_{i,1}\otimes v_{i,2}\otimes v_{i,3}=0$ and the
right-hand side of \eqref{eq:segre-differential} is
$V_1\otimes V_2\otimes V_3$, then
\[
 \MaxR_\F(\dim V_1,\dim V_2,\dim V_3)\le r.
\]
\end{corollary}
\begin{proof}
Apply Lemma~\ref{lem:regular-zero} to the affine parameter space $Y$.
\end{proof}

\section{Cancellation and amplification}\label{sec:amplification}

To turn a regular tuple into a regular cancellation,
we tensor it with an auxiliary cancellation. The cancellation must
preserve independent variations in all three factors. The following
lemma provides exactly this property.

\Needspace{8\baselineskip}
\begin{lemma}[Spanning cancellation]\label{lem:cancellation}
Let $\F$ be a field with more than two elements and $k\ge1$. Then there are $k^2+k$ triples
$(x_\nu,y_\nu,z_\nu)$ in three $k$-dimensional spaces such that
\begin{equation}\label{eq:cancellation}
 \sum_{\nu=1}^{k^2+k}x_\nu\otimes y_\nu\otimes z_\nu=0,
\end{equation}
and each of the families $(x_\nu\otimes y_\nu)_\nu$,
$(x_\nu\otimes z_\nu)_\nu$, and $(y_\nu\otimes z_\nu)_\nu$
spans its two-factor tensor product.
\end{lemma}
\begin{proof}
Let $\mathcal{A}=\F^k$ with coordinatewise multiplication. Let $e_\ell$ and $\epsilon_\ell$ be the
coordinate bases of $\mathcal{A}$ and $\mathcal{A}^*$. Fix two distinct nonzero elements $a,b\in\F$. For $1\leq i\leq k$ let $v_i:=b e_i+\sum_{j\neq i}ae_j$ and $v_0:=a\sum_{j=1}^k e_j$. We note that $v_i-v_0=(b-a)e_i$ and thus $v_0,\dots,v_k$ linearly span $\mathcal{A}$. In particular there exists a basis
$u_1,\ldots,u_k$ of $\mathcal{A}$ consisting of units, i.e.~each basis vector has all coordinates nonzero. Let $u_1^*,\ldots,u_k^*$
be the dual basis. The multiplication tensor has the two expressions
\begin{equation}\label{eq:multiplication}
 \sum_{i,j=1}^k u_i^*\otimes u_j^*\otimes u_i u_j
 =\sum_{\ell=1}^k \epsilon_\ell\otimes\epsilon_\ell\otimes e_\ell.
\end{equation}
Subtracting gives \eqref{eq:cancellation} in $ \mathcal{A}^*\otimes \mathcal{A}^*\otimes \mathcal{A}$.

The first $k^2$ summands already have all three spanning properties.
The tensors $u_i^*\otimes u_j^*$ form a basis of $\mathcal{A}^*\otimes \mathcal{A}^*$.
For fixed $j$, the vectors $(u_i u_j)_i$ form a basis of $\mathcal{A}$,
because multiplication by the unit $u_j$ is invertible. Hence
$(u_i u_j)\otimes u_j^*$ span $\mathcal{A}\otimes \mathcal{A}^*$; the other pair is
handled in the same way.
\end{proof}
A cancellation satisfying the conclusions of Lemma \ref{lem:cancellation} is called \emph{spanning}.

%The identity uses only a split algebra and a basis of units, and
%works in every characteristic. For related uses of multiplication
%tensors in tensor geometry, see \cite{LandsbergMichalekAbelian}.

We recall the parameterisation of tensors in format $n_1\times n_2\times n_3$ of rank at most $r$:
\begin{equation}\label{eq:map-r} 
\Phi_r: (\F^{n_1}\times \F^{n_2}\times\F^{n_3})^r\rightarrow \F^{n_1}\otimes \F^{n_2}\otimes\F^{n_3}.
\end{equation}
\begin{proposition}[Amplification]\label{prop:amplification}
Let $\F=\R$, or let $\F$ be algebraically closed. Suppose that for some $r,n_1,n_2,n_3$ the map \eqref{eq:map-r}
has surjective differential at some point. Then, for every $k\ge1$,
\begin{equation}\label{eq:amplification}
 \MaxR_\F(kn_1,kn_2,kn_3)\le r(k^2+k).
\end{equation}
\end{proposition}
\begin{proof}
Let $(v_{\mu,j})_\mu$ be the given regular point, with $v_{\mu,j}\in \F^{n_j}$,
and choose a spanning cancellation $(x_{\nu},y_{\nu},z_{\nu})_\nu$ as in
Lemma~\ref{lem:cancellation}, with $x_{\nu},y_\nu,z_\nu\in \F^k$. In the enlarged spaces $(\F^{n_j}\otimes \F^k)$, take the
$r(k^2+k)$ rank one tensors:
\[
 (v_{\mu,1}\otimes x_{\nu})\otimes
 (v_{\mu,2}\otimes y_{\nu})\otimes
 (v_{\mu,3}\otimes z_{\nu}).
\]
We will prove that they form a regular cancellation.
Indeed, their sum is zero, since for every fixed $\mu$, after regrouping, we have:
\[\sum_\nu (v_{\mu,1}\otimes x_{\nu})\otimes (v_{\mu,2}\otimes y_{\nu})\otimes (v_{\mu,3}\otimes z_{\nu})=v_{\mu,1}\otimes v_{\mu,2}\otimes v_{\mu,3}\otimes\left(\sum_\nu x_{\nu}\otimes  y_{\nu}\otimes z_{\nu}\right)=0. \]

The spanning properties of Lemma \ref{lem:cancellation} ensure that each original
tangent direction is tensored with the entire auxiliary space.
Indeed, fixing $\mu$ and varying the first enlarged factor, independently
for every $\nu$, gives
\[
 (\F^{n_1}\otimes v_{\mu,2}\otimes v_{\mu,3})
       \otimes(\F^k\otimes \F^k\otimes \F^k),
\]
under the natural reordering of factors. Varying the other two factors
gives the analogous subspaces. By
\eqref{eq:segre-differential}, their sum over $\mu$ is the entire enlarged
tensor space. The constructed cancellation in $(\F^{kn_1}\times \F^{kn_2}\times \F^{kn_3})$ is therefore regular, and
Corollary~\ref{cor:segre-zero} proves \eqref{eq:amplification}.
\end{proof}

\section{The asymptotic tensor bound}\label{sec:segre-proof}

To finish the proof of the main result, we need a good understanding of regular tuples.
We use the tangent-space formulation of Lickteig's construction;
this formulation matters in positive characteristic.

\begin{lemma}[Lickteig's estimate]\label{lem:lickteig}
Let $\F$ be infinite, and fix $a,b,c>0$. For all sufficiently large $s$,
the decomposition map in format
$\lceil as\rceil\times\lceil bs\rceil\times\lceil cs\rceil$
has an $\F$-point with surjective differential using
\[
 r_s=\frac{abc}{a+b+c}s^2+O_{a,b,c}(s)
\]
summands.
\end{lemma}
\begin{proof}
The finite estimate proved in Lemma~\ref{lem:rectangular-tangents}
is the following: for $2\le p_1\le p_2\le p_3$, with $p_2,p_3$ even,
there is a surjective differential with
\[
 r=p_3\left\lceil\frac{p_1p_2}{p_1+p_2+p_3-2}\right\rceil
\]
summands. This is Lickteig's rectangular estimate in the precise
form needed here.

Order the dimensions $\lceil as\rceil,\lceil bs\rceil,\lceil cs\rceil$,
and round the last two upwards to even integers. The enlarged
dimensions satisfy $p_1=as+O(1)$, $p_2=bs+O(1)$, $p_3=cs+O(1)$,
after the same permutation of $a,b,c$. Hence
\[
 r=\frac{p_1p_2p_3}{p_1+p_2+p_3-2}+O(p_3)
  =\frac{abc}{a+b+c}s^2+O_{a,b,c}(s).
\]

Fix linear projections $\F^{p_1}\rightarrow \F^{\lceil as\rceil}$, $\F^{p_2}\rightarrow \F^{\lceil bs\rceil}$, $\F^{p_3}\rightarrow \F^{\lceil cs\rceil}$ and induce from them surjective linear maps:
%\jbI{tu powyżej trochę przeredagowałem zdanie. Jeszcze: stricte, to poniżej to nie są liniowe odwzorowania, bo $\times$ nie jest operacją na przestrzeniach liniowych, tu bierzemy produkt rozmaitości, ale może nie trzeba tego wyjaśniać}:
\begin{align*}P:(\F^{p_1}\times \F^{p_2}\times\F^{p_3})^r\rightarrow (\F^{\lceil as\rceil}\times \F^{\lceil bs\rceil}\times\F^{\lceil cs\rceil})^r, \text{ and }\\
\pi:\F^{p_1}\otimes\F^{p_2}\otimes \F^{p_3}\rightarrow \F^{\lceil as\rceil}\otimes \F^{\lceil bs\rceil}\otimes\F^{\lceil cs\rceil}.\end{align*}
Let
\begin{align*}
\Phi_r:(\F^{p_1}\times \F^{p_2}\times\F^{p_3})^r\rightarrow \F^{p_1}\otimes\F^{p_2}\otimes \F^{p_3}\text{ and }\\
\widetilde\Phi_r:(\F^{\lceil as\rceil}\times \F^{\lceil bs\rceil}\times\F^{\lceil cs\rceil})^r\rightarrow \F^{\lceil as\rceil}\otimes \F^{\lceil bs\rceil}\otimes\F^{\lceil cs\rceil}
\end{align*}
be the two maps parameterising tensors of rank at most $r$.
We have:
\[
 \widetilde\Phi_r\circ P=\pi\circ \Phi_r.
\]
At a regular tuple, the differential of the right-hand
side is surjective. The chain rule therefore makes the differential
of $\widetilde\Phi_r$ surjective at the projected tuple. This proves the claim.
\end{proof}

We emphasize that the above argument proves a surjective differential over the field of any characteristic.
It does not deduce this from dominance of the maps $\Phi_r$ or $\widetilde{\Phi}_r$,
as such an implication can
fail in positive characteristic.%\jb{tu zmieniłem}

\begin{proof}[Proof of Theorem~\ref{thm:segre}]
Let $\gamma=abc/(a+b+c)$,
$s=\lfloor\sqrt n\rfloor$, and $k=\lceil n/s\rceil$.
Lemma~\ref{lem:lickteig} supplies a regular $r_s$-tuple of the format $\lceil as\rceil\times\lceil bs\rceil\times\lceil cs\rceil$%\jb{tu dodałem}
 with
$r_s=\gamma s^2+O_{a,b,c}(s)$. Since $sk\ge n$, each
factor of the desired format is a quotient of the corresponding
factor of the amplified format. Tensor rank does not increase under a
linear map on the factors. Proposition~\ref{prop:amplification}
therefore gives
\[
 \MaxR_\F(\lfloor an\rfloor,\lfloor bn\rfloor,\lfloor cn\rfloor)
 \le(k^2+k)r_s
 =\gamma n^2+O_{a,b,c}(n^{3/2}),
\]
where $sk=n+O(\sqrt n)$.

For the lower bound, the cone of decomposable tensors in format
$n_1\times n_2\times n_3$ has dimension $n_1+n_2+n_3-2$.
The closure of its $r$-fold sum has dimension at most
$r(n_1+n_2+n_3-2)$. Consequently,
\begin{equation}\label{eq:dimension-lower}
 \MaxR_\F(n_1,n_2,n_3)\ge
 \left\lceil\frac{n_1n_2n_3}{n_1+n_2+n_3-2}\right\rceil.
\end{equation}
%For $\R$, the same conclusion follows because a proper complex
%algebraic subset cannot contain every real tensor.
 Substituting
$n_1=\lfloor an\rfloor$, $n_2=\lfloor bn\rfloor$,
$n_3=\lfloor cn\rfloor$ gives $\gamma n^2+O_{a,b,c}(n)$,
and completes the proof.
\end{proof}

\section{An explicit bound for equal dimensions}\label{sec:effective}

%For $n\times n\times n$ tensors, Lickteig's exact tangent-space count gives a
%finite bound. Optimizing it requires some care: rounding two
%independently chosen parameters of order $\sqrt n$ can lose the
%second-order coefficient.
In this section we provide an explicit upper bound on the error term and hence on maximal rank.

\Needspace{16\baselineskip}
\begin{lemma}[Effective bound for equal dimensions]\label{lem:effective}
Let $\F=\R$, or let $\F$ be algebraically closed. For every $n\ge64$,
\begin{equation}\label{eq:closed-bound}
 \MaxR_\F(n,n,n)\le
 \frac{n^2}{3}+\frac{2\sqrt6}{9}n^{3/2}+\frac83n^{5/4}+\frac{65}{12}n.
\end{equation}
A sharper finite bound, valid for every $n\ge1$, is
\begin{equation}\label{eq:finite-bound}
 \MaxR_\F(n,n,n)\le B(n),\qquad
 B(n)=\min_{1\le s\le n}g(s)
       \left\lceil\frac ns\right\rceil
       \left(\left\lceil\frac ns\right\rceil+1\right),
\end{equation}
where $g(3)=5$ and $g(s)=\lceil s^3/(3s-2)\rceil$ for $s\ne3$.
\end{lemma}
The coefficient $2\sqrt6/9$ is optimal for the function $B$. We prove below that
\begin{equation}\label{eq:B-lower}
 B(n)\ge\frac{n^2}{3}+\frac{2\sqrt6}{9}n^{3/2}+\frac{2n}{9}.
\end{equation}
\begin{proof}[Proof of Lemma~\ref{lem:effective}]
By Lemma~\ref{lem:equal-tangents} and
Proposition~\ref{prop:amplification}, followed by projection from
dimension $sk$ to dimension $n$ in each factor,
\[
 \MaxR_\F(n,n,n)\le g(s)k(k+1)
 \qquad\text{for positive integers }s,k\text{ with }sk\ge n.
\]
For fixed $s$, the smallest admissible $k$ is $\lceil n/s\rceil$.
Minimizing over $1\le s\le n$ proves \eqref{eq:finite-bound}.

We first estimate $g(s)$. Polynomial division gives
\[
 \frac{s^3}{3s-2}
 =\frac{s^2}{3}+\frac{2s}{9}+\frac{4s}{9(3s-2)}.
\]
For $s\ge2$, the inequality $4s<3(3s-2)$ shows that the last
term lies in $(0,1/3)$. Applying $x\le\lceil x\rceil<x+1$
for $s\ne3$, and checking $g(1)=1$ and $g(3)=5$ separately, gives
\begin{equation}\label{ineq1}
 \frac{s^2}{3}+\frac{2s}{9}<g(s)
 \le\frac{s^2}{3}+\frac{2s}{9}+\frac43
 \qquad(s\ge1).
\end{equation}
For positive integers $s,k$, put $N=sk$ and define
\begin{equation}\label{eq:effective-E}
 E(s,k)=\left(\frac{s^2}{3}+\frac{2s}{9}\right)k(k+1)
       =\frac{N^2}{3}+\frac{N(3s+2k)}9+\frac{2N}{9}.
\end{equation}
Multiplying \eqref{ineq1} by $k(k+1)$ yields
\begin{equation}\label{eq:E-bounds}
 E(s,k)\le g(s)k(k+1)\le E(s,k)+\frac43(k^2+k).
\end{equation}
The AM--GM inequality applied to $3s$ and $2k$ gives
$3s+2k\ge2\sqrt{6N}$. Substituting this into
\eqref{eq:effective-E} and using the lower inequality in
\eqref{eq:E-bounds}, we obtain the lower bound
\eqref{eq:B-lower} by replacing $N\ge n$ with $n$ and minimizing
over the pairs $s,k$ defining $B(n)$.

For the upper bound, assume $n\ge64$. We will choose $s,k$ by first
fixing $t=3s+2k$ close to $2\sqrt{6n}$, and then rounding a root
of $x(t-3x)/2=n$. Set
\[
 t=\left\lceil\sqrt{24n+36}\right\rceil,
 \qquad d=\sqrt{t^2-24n}.
\]
Since $n\ge64>54$, it follows that $2\sqrt{6n}>36$.
Squaring therefore shows that
$\sqrt{24n+36}<2\sqrt{6n}+1/2$.
Combining this with $\lceil x\rceil<x+1$ gives
\begin{equation}\label{eq:effective-t}
 t<2\sqrt6\,\sqrt n+\frac32.
\end{equation}
The definition of $t$ gives $d^2\ge36$.
For an upper bound, squaring \eqref{eq:effective-t} and subtracting $24n$
gives $d^2<6\sqrt6\,\sqrt n+9/4$.
Using $\sqrt6<5/2$ and $\sqrt n\ge8$, we can bound the two summands by $15\sqrt n$ and $(9/32)\sqrt n$, respectively.
Since $15+9/32<16$, taking square roots yields
\begin{equation}\label{eq:effective-d}
 6\le d<4n^{1/4}.
\end{equation}
Since $t^2-d^2=24n>0$, we also have $d<t$.

The roots of $x(t-3x)/2=n$ are
\[
 \alpha=\frac{t-d}{6},\qquad \beta=\frac{t+d}{6}.
\]
By \eqref{eq:effective-d}, their distance $d/3$ is at least $2$.
Choose the largest integer $s\equiv t\pmod2$ with $s\le\beta$.
Consecutive integers of this parity differ by $2$, so
$\delta:=\beta-s$ satisfies $0\le\delta<2$, and consequently
$\alpha\le s\le\beta$.
Since $d<t$, both roots lie in $(0,t/3)$.
Define $k=(t-3s)/2$; the parity condition makes $k$ an integer,
and $s<t/3$ makes it positive.
To check $s\le n$, note that $n\ge64$ gives
$\sqrt n\le n/8$ and $3/2\le n/8$. Together with $2\sqrt6<5$,
these bounds turn \eqref{eq:effective-t} into $t<3n/4$,
so indeed $s<t/3<n$.

Substituting $s=\beta-\delta$ into $N=s(t-3s)/2$ gives
\[
 N-n=\frac{\delta(d-3\delta)}2.
\]
Because $0\le\delta<2$ and $d\ge6$, the right-hand side is
nonnegative and is at most $\delta d/2<d$.
Together with \eqref{eq:effective-d}, this proves
\begin{equation}\label{eq:effective-N}
 n\le N<n+4n^{1/4}.
\end{equation}

The same substitution gives $k=(t-d)/4+3\delta/2$.
Using $d\ge6$ and $\delta<2$ first, and then
\eqref{eq:effective-t}, bounds $k$ by
\[
 k<\frac t4+\frac32
   <\frac{\sqrt6}{2}\sqrt n+\frac{15}{8}.
\]
Since $\sqrt6<5/2$ and $\sqrt n\ge8$, the two terms in the
final expression are smaller than $(5/4)\sqrt n$ and
$(1/4)\sqrt n$, respectively. Adding these bounds gives
\begin{equation}\label{eq:effective-k}
 k<\frac32\sqrt n.
\end{equation}

It remains to estimate $N^2$, $Nt$, $2N$, and $k^2+k$, which
occur in \eqref{eq:E-bounds} after expanding $E(s,k)$ and using
$3s+2k=t$.
First, square \eqref{eq:effective-N} and expand.
Since $\sqrt n\ge8$ implies $16\sqrt n\le2n$, this gives
\begin{equation}\label{eq:effective-N2}
 N^2\le n^2+8n^{5/4}+2n.
\end{equation}
Next, $n\ge64$ implies $n^{1/4}\ge\sqrt6$ and $n^{3/4}\ge8$,
so the two cross terms in the product of
\eqref{eq:effective-N} and \eqref{eq:effective-t} satisfy
\[
 8\sqrt6\,n^{3/4}\le8n,\qquad 6n^{1/4}\le n.
\]
Multiplying \eqref{eq:effective-N} and \eqref{eq:effective-t},
expanding, and applying these cross-term estimates gives
\begin{equation}\label{eq:effective-Nt}
 Nt\le2\sqrt6\,n^{3/2}+\left(\frac32+8+1\right)n
     =2\sqrt6\,n^{3/2}+\frac{21}{2}n.
\end{equation}
Similarly, $n^{3/4}\ge8$ implies $8n^{1/4}\le n$.
Doubling \eqref{eq:effective-N} consequently yields
\begin{equation}\label{eq:effective-2N}
 2N\le3n.
\end{equation}
Finally, \eqref{eq:effective-k} gives
$k^2+k<(9/4)n+(3/2)\sqrt n$.
Using $\sqrt n\le n/8$ in the second term gives
\begin{equation}\label{eq:effective-kk}
 k^2+k\le\frac{39}{16}n.
\end{equation}

The chosen pair satisfies $1\le s\le n$ and $sk\ge n$, so
$k\ge\lceil n/s\rceil$ and the definition of $B(n)$ gives
$B(n)\le g(s)k(k+1)$.
Now use \eqref{eq:effective-E} with $3s+2k=t$ in
\eqref{eq:E-bounds}, and apply
\eqref{eq:effective-N2}--\eqref{eq:effective-kk}. This gives
\[
 \begin{aligned}
 B(n)&\le\frac{N^2}{3}+\frac{Nt}{9}+\frac{2N}{9}
          +\frac43(k^2+k)\\
 &\le\frac{n^2}{3}+\frac{2\sqrt6}{9}n^{3/2}
       +\frac83n^{5/4}
       +\left(\frac23+\frac76+\frac13+\frac{13}{4}\right)n\\
 &=\frac{n^2}{3}+\frac{2\sqrt6}{9}n^{3/2}
       +\frac83n^{5/4}+\frac{65}{12}n.
 \end{aligned}
\]
Since $n^{5/4}=o(n^{3/2})$, this upper bound and
\eqref{eq:B-lower} also prove the asserted optimality of
$2\sqrt6/9$ for $B(n)$.
\end{proof}

\appendix
\section{The tangent-space estimates}\label{app:tangents}

For a nonzero decomposable tensor $t=u\otimes v\otimes w\in A\otimes B\otimes C$,
write
\[
 \mathcal T_{A,B,C}(t)=\mathcal T(t)
 =A\otimes v\otimes w+u\otimes B\otimes w+u\otimes v\otimes C.
\]
This is the tangent space to the affine Segre cone at $t$.
By \eqref{eq:segre-differential}, a tuple gives a surjective
differential precisely when its tangent spaces span the ambient space.

We first explain the passage from an algebraic closure to an infinite
ground field $\F$. For fixed $A,B,C$ and $r$, consider
\[
 P=\bigl((A\setminus\{0\})\times(B\setminus\{0\})
       \times(C\setminus\{0\})\bigr)^r.
\]
The condition that the corresponding tangent spaces span $A\otimes B\otimes C$
defines an open subset $U\subseteq P$ over $\F$: it is a full-rank
condition on a matrix whose entries are polynomial functions of the
factors. If a spanning tuple exists over $\overline{\F}$, then the product of minors of this matrix must be a nonzero polynomial. Since $\F$ is infinite and $P$ is open in an affine space, this product of minors cannot vanish on all $F$-points of $P$. Hence, if a spanning tuple exists over $\overline{\F}$, then $U$ is nonempty. Thus it suffices to prove the two existence
statements below over $\overline{\F}$.

\subsection*{The rectangular estimate}

\begin{lemma}[Rectangular tangent span]
\label{lem:rectangular-tangents}
Let $\F$ be infinite, let $2\le a\le b\le c$, and suppose that $b,c$
are even. In format $a\times b\times c$, there are
\begin{equation}\label{eq:rectangular-count}
 r=c\left\lceil\frac{ab}{a+b+c-2}\right\rceil
\end{equation}
nonzero decomposable tensors whose tangent spaces to the affine Segre
cone span the ambient space.
\end{lemma}
\begin{proof}
Over an algebraically closed field of arbitrary characteristic,
the proof of \cite[Thm~3.3, pp.~105--107]{Lickteig} constructs
$r$ tensors with the required tangent-space property, assuming
$a\le b\le c$ and $b,c$ even. The preceding argument gives such a
tuple over $\F$.
\end{proof}

\subsection*{The sharp count in equal dimensions}

\begin{lemma}[Lickteig; equal-dimensional tangent span]
\label{lem:equal-tangents}
Let $\F$ be infinite and $s\ge1$. Set
\[
 g(3)=5,\qquad
 g(s)=\left\lceil\frac{s^3}{3s-2}\right\rceil\quad(s\ne3).
\]
In format $s\times s\times s$, there are $g(s)$ nonzero decomposable
tensors whose tangent spaces to the affine Segre cone span the ambient space.
\end{lemma}
\begin{proof}
By the argument above, we may assume that $\F$ is algebraically closed.
Let $V=A\otimes B\otimes C$, where $\dim A=\dim B=\dim C=s$.
For $s\ge5$, $s\ne7$, the construction in the proof of
\cite[Thm~4.4, pp.~110--112]{Lickteig} gives $g(s)$ tensors
whose tangent spaces span $V$, in every characteristic.

For $s=7$, take a tensor-product subspace $M\subset V$ of format
$6\times6\times5$. \cite[Prop.~4.2]{Lickteig}, with $n=6$,
and Lemma~4.7, with $p=6$ and $q=18$, give twelve tensors in $M$
and six further tensors in $V$ whose tangent sum $H$ contains $M$ and satisfies
\[
 \dim H=12\cdot15+4\cdot12+19\cdot6=342.
\]
Since $\dim V=343$ and decomposable tensors span $V$, choose a
decomposable tensor $t\notin H$. Then $t\in\mathcal T_{A,B,C}(t)$, so
$H+\mathcal T_{A,B,C}(t)=V$,%\jb{tu zmieniłem, żeby było spójnie z oznaczeniem na początku sekcji},
 using $19=g(7)$ tensors.

For $s=1,2$, the $s$ diagonal tensors $e_i\otimes e_i\otimes e_i$
already have spanning tangent spaces, and $g(s)=s$.
The cases $s=3,4$ are proved explicitly in the following subsection.
\end{proof}

\subsection*{The small-dimensional cases}

We give the constructions used above for $s=3,4$. In each factor,
choose a basis $e_1,\ldots,e_s$. The sum of the tangent spaces at
\[
 e_i\otimes e_i\otimes e_i,\qquad 1\le i\le s,
\]
is the subspace spanned by the coordinate tensors with at least two
equal indices. Indeed, varying one factor at the $i$th diagonal tensor
gives precisely the coordinate tensors in which the other two indices
are equal to $i$. The quotient by this subspace therefore has basis
\[
 e_i\otimes e_j\otimes e_k,
 \qquad i,j,k\text{ pairwise distinct}.
\]

For $s=3$, add
\[
 e_1\otimes e_2\otimes e_3,
 \qquad u\otimes u\otimes u,
 \qquad u=e_1+e_2+e_3.
\]
In the six-dimensional quotient just described, the tangent space at
$e_1\otimes e_2\otimes e_3$ has image equal to the line spanned by
this tensor: every other coordinate variation has a repeated index.

After quotienting by this line, order the remaining coordinates as
\[
 (1,3,2),\ (2,1,3),\ (2,3,1),\ (3,1,2),\ (3,2,1).
\]
The five vectors
\[
 e_1\otimes u\otimes u,\quad
 e_2\otimes u\otimes u,\quad
 e_3\otimes u\otimes u,\quad
 u\otimes e_1\otimes u,\quad
 u\otimes e_2\otimes u
\]
belong to the tangent space at $u\otimes u\otimes u$. Their
coordinate columns in this order form the matrix
\[
 \begin{pmatrix}
  1&0&0&0&0\\
  0&1&0&1&0\\
  0&1&0&0&0\\
  0&0&1&1&0\\
  0&0&1&0&1
 \end{pmatrix}.
\]
Its determinant is $1$, so these vectors span the remaining quotient
over every field. The tangent spaces at the five chosen tensors
therefore span the whole space, as required for $g(3)=5$.

For $s=4$, put $u=e_1+e_2+e_3+e_4$ and $e_{ij}=e_i+e_j$.
In addition to the four diagonal tensors, take
\[
 u\otimes u\otimes u,\qquad
 p=e_{12}\otimes e_{23}\otimes e_{34},\qquad
 q=e_{13}\otimes e_{24}\otimes e_{14}.
\]
For $t=a\otimes b\otimes c$, write
\[
 \alpha_i(t)=e_i\otimes b\otimes c,\qquad
 \beta_i(t)=a\otimes e_i\otimes c,\qquad
 \gamma_i(t)=a\otimes b\otimes e_i.
\]
These vectors span $\mathcal T(t)$.

Let $Q$ be the quotient by the sum of the tangent spaces at the
four diagonal tensors and at $p,q$, and denote the class of
$e_i\otimes e_j\otimes e_k$ by $\tau_{ijk}$. Classes with repeated
indices vanish. Expanding the coordinate variations at $p$ gives
\[
 \tau_{123}=\tau_{124}=\tau_{132}=\tau_{134}
 =\tau_{231}=\tau_{234}=\tau_{324}=\tau_{423}=0,
\]
together with
\[
 \tau_{213}+\tau_{214}=0,\qquad
 \tau_{143}+\tau_{243}=0.
\]
The variations $\alpha_2,\alpha_4,\beta_1,\beta_4,\beta_2,
\gamma_3$ at $q$ then give, respectively,
\[
 \tau_{241}=\tau_{421}=\tau_{314}=\tau_{341}
 =\tau_{321}=\tau_{143}=0,
\]
and $\gamma_2(q)$ gives $\tau_{142}+\tau_{342}=0$.
Hence $\tau_{243}=0$ as well, and $Q$ is generated by
\[
 (v_1,\ldots,v_7)
 =(\tau_{214},\tau_{142},\tau_{431},\tau_{432},
   \tau_{312},\tau_{412},\tau_{413}).
\]
The following variations at $u\otimes u\otimes u$ have the
indicated images in $Q$:
\[
\begin{array}{c|ccccccc}
 \text{variation}&\gamma_4&\alpha_1&\gamma_1&\beta_3
    &\alpha_3&\gamma_2&\gamma_3\\ \hline
 \text{image}&v_1&v_2&v_3&v_3+v_4
    &v_5-v_2&v_4+v_5+v_6&v_7-v_1
\end{array}
\]
Read from left to right, these images successively span
$v_1,\ldots,v_7$. The seven chosen tangent spaces therefore span
the ambient space. Since $g(4)=7$ and the computation uses only
addition and subtraction, this proves the assertion over every field.

\subsection*{Acknowledgements}
Buczy{\'n}ski is supported by ``Advanced problems in contact manifolds, secant varieties, and their generalisations (\mbox{APRICOTS+})'', the project 2023/51/B/ST1/02799 awarded by National Science Center, Poland. . Micha{\l}ek is supported by the DFG grant 580118961.

\subsection*{Use of artificial intelligence}
Throughout our work we got a great help from ChatGPT Astra, which in particular provided the key ideas, together with their exposition. The article arose after many discussions of the authors, also with other colleagues. We believe the ideas presented here are very much worth sharing.

\end{document}